  \documentstyle{amsppt}


  \loadmsbm
  \loadbold


\define\dQ{{\Bbb Q}}

\define\dZ{{\Bbb Z}}

\define\cO{{\Cal O}}

\define\ep{{\frak p}}

\define\eq{{\frak q}}


\define\Gal{{\operatorname{Gal}}}

\define\Pic{{\operatorname{Pic}}}

\define\tw#1{\,^{#1}\!}

\font\pfeile = cmsy10 scaled 1440
\newfam\pfeilfam
\textfont\pfeilfam=\pfeile
                   \scriptfont\pfeilfam=\pfeile
                                      \scriptscriptfont\pfeilfam=\pfeile

\mathchardef\swpfeil="2D2E

\def\\{\let\stoken= } \\
\long\def\unexpandedwrite#1#2{\def\finwrite{\write#1}%
{\aftergroup\finwrite\aftergroup{\sanitize#2\endsanity}}}
\def\sanitize{\futurelet\next\sanswitch}
\def\sanswitch{\ifx\next\endsanity
\else\ifcat\noexpand\next\stoken\aftergroup\space\let\next=\eat
\else\ifcat\noexpand\next\bgroup\aftergroup{\let\next=\eat
\else\ifcat\noexpand\next\egroup\aftergroup}\let\next=\eat
\else\let\next=\copytoken\fi\fi\fi\fi \next}
\def\eat{\afterassignment\sanitize \let\next= }
\long\def\copytoken#1{\ifcat\noexpand#1\relax\aftergroup\noexpand
\else\ifcat\noexpand#1\noexpand~\aftergroup\noexpand\fi\fi
\aftergroup#1\sanitize}
\def\endsanity\endsanity{}

\expandafter\ifx\csname pre numero.tex at\endcsname\relax \else
\endinput\fi \expandafter\chardef\csname pre numero.tex
at\endcsname=\the\catcode`\@ \catcode`\@=11
\def\tokenitize{\futurelet\next\tokenswitch}
\def\tokenswitch{\ifx\next\endtokenity
\else\ifcat\noexpand\next\stoken\aftergroup\space\let\next=\eatt
\else\let\next=\copytok\fi\fi\next}
\def\eatt{\afterassignment\tokenitize \let\next= }
\long\def\copytok#1{\aftergroup\string\aftergroup#1\tokenitize}
\def\endtokenity{}
\newcount \secno
\newcount \Refno
\def\phantomlabel#1{{\aftergroup\expandafter\aftergroup\xdef%
\aftergroup\csname\tokenitize ref#1\endtokenity\aftergroup\endcsname}%
{{\the\secno.\the\Refno}}%
\global\advance\Refno by 1\relax}
\def\label#1{{\the\secno.\the\Refno}\phantomlabel{#1}%
\unexpandedwrite\labx{\phantomlabel{#1}}}
\def\Ref#1{{\aftergroup\expandafter\aftergroup\ifx
\aftergroup\csname\tokenitize ref#1\endtokenity}\endcsname\relax
??R\'ef\'erence {\tokenitize#1\endtokenity} non d\'efinie??%
\message{Reference #1 non definie}%
\else{\aftergroup\csname\tokenitize ref#1\endtokenity}\endcsname\fi
\ifnextsecn@\else\message{ATTENTION secno n'a pas ete initialise}\fi }
\def\phantomsection{\global\advance\secno by 1 \global\Refno = 1\relax}
\def\secnum{{\phantomsection\the\secno\unexpandedwrite\labx{\phantomsection}}}
\newif\ifnextsecn@
\def\nextsecno#1{\nextsecn@true\global\secno=#1\global\advance\secno by -1%
\immediate\write\labx{\secno=\the\secno}\message{\the\secno} }
\openin1 \jobname.lab \ifeof1 \message
                                {No file \jobname.lab.}\else\closein1\relax\input \jobname.lab \fi
\newwrite\labx
\immediate\openout\labx=\jobname.lab 
\catcode`\@=\csname pre numero.tex at\endcsname


\scrollmode \NoBlackBoxes \magnification=\magstep1 \nologo
        \pagewidth{6.5truein}
        \pageheight{9truein}\topmatter

\parskip=5pt


\TagsOnRight
\define\albert{1}

\define\as{2}
\define\AT{3}
\define\FS{4}
\define\FKS{5}
\define\FSi{6}
\define\hayes{7}
\define\KS{8}
\define\Po{9}
\define\Weil{10}

\define\charr{\text{\rm{char}}}

\define\ov{\overline}

\define\wt{\widetilde}
\define\el{\frak l}
\define\ea{\frak a}
\define\eb{\frak b}

\title
Abelian extensions of global fields with constant local degrees \endtitle

\author
Hershy Kisilevsky and Jack Sonn \endauthor \affil Concordia University,
Montreal and Technion, Haifa \endaffil
\address
Department of Mathematics, Concordia University, 1455 de Maisonneuve
Blvd. West, Montreal, Canada, and Department of Mathematics, Technion,
32000 Haifa, Israel \endaddress \email kisilev\@mathstat.concordia.ca,
sonn\@math.technion.ac.il \endemail
\thanks
The research of the first author was supported by NSERC and FQRNT.
The research of the second author was supported by the VPR
Fund at the Technion.
\endthanks

\abstract Given a global field $K$ and a positive integer $n$, there
exists an abelian extension $L/K$ (of exponent $n$) such that the local
degree of $L/K$ is equal to $n$ at every finite prime of $K$, and is
equal to two at the real primes if $n=2$. As a consequence, the
$n$-torsion subgroup of the Brauer group of $K$ is equal to the relative
Brauer group of $L/K$. \endabstract

\endtopmatter
\document

\head 1. Introduction \endhead

Let $K$ be a field, $Br(K)$ its Brauer group.  If $L/K$ is a field
extension, then the relative Brauer group $Br(L/K)$ is the kernel of the
restriction map $res_{L/K}:Br(K)\rightarrow Br(L)$. Relative Brauer
groups have been studied by Fein and Schacher (see e.g.
\cite{\FS,\FKS,\FSi}.)  Every subgroup of $Br(K)$ is a relative Brauer
group $Br(L/K)$ for some extension $L/K$ \cite{\FS}, and the question
arises as to which subgroups of $Br(K)$ are \it algebraic relative Brauer
groups,  i.e. \rm of the form $Br(L/K)$ with $L/K$ an algebraic
extension. For example if $L/K$ is a finite extension of number fields,
then $Br(L/K)$ is infinite \cite{\FKS}, so no finite subgroup of $Br(K)$
is an algebraic relative Brauer group. In \cite \as \ the question was
raised as to whether or not the $n$-torsion subgroup $Br_n(K)$ of the
Brauer group $Br(K)$ of a field $K$ is an  algebraic relative Brauer
group.  For example, if $K$ is a ($p$-adic) local field, then
$Br(K)\cong\dQ/\dZ$, so $Br_n(K)$ is an algebraic relative Brauer group
for all $n$.
  A counterexample was given in
\cite \as \ for $n=2$ and $K$ a formal power series field over a local
field.   For global fields $K$, the problem is a purely arithmetic one,
because of the fundamental local-global description of the Brauer group
of a number field. In particular, for a Galois extension $L/K$ of global
fields, if the local degree of $L/K$ at every finite prime is equal to
$n$, and is equal to 2 at the real primes for $n$ even, then
$Br(L/K)=Br_n(K)$.   In \cite \as, it was proved that $Br_n(\dQ)$ is an
algebraic relative Brauer group for all squarefree $n$. In \cite \KS,
 the arithmetic criterion above was verified for
 any number field $K$ Galois over $\dQ$ and any $n$ prime to the
class number of $K$, so in particular, $Br_n(\dQ)$ is an algebraic
relative Brauer group for all $n$.  In \cite \Po, Popescu proved that for
a global function field $K$ of characteristic $p$, the arithmetic
criterion holds for $n$ prime to the order of the non-$p$ part of the
Picard group of $K$. \smallskip In this paper we settle the question
completely, by verifying the arithmetic criterion for all $n$ and all
global fields $K$. In particular,  the $n$-torsion subgroup of the Brauer
group of $K$ is an algebraic relative Brauer group for all $n$ and all
global fields $K$. The proof, an extension of the ideas in \cite \KS,
reduces to the case $n$ a prime power $\ell^r$.  We first carry out the
proof for number fields $K$.  The proof for the function field case when
$\ell\neq \charr(K)$ is essentially the same as the proof in the number
field case. The proof for $\ell=\charr(K)$ appears in \cite \Po.

\head 2. A Splitting Lemma
    \endhead
 Let $K$ be a number field, $\ep$ a finite prime of $K$, $I_{\ep}$ the
group of fractional
 ideals prime to $\ep$, $P_{\ep}$ the group of principal fractional
ideals in $I_{\ep}$,
 $P_{\ep,1}$ the group of principal fractional ideals in $P_{\ep}$
generated by elements
 congruent to $1$ mod $\ep$.  Then $Cl_K\cong I_{\ep}/P_{\ep}$ is the
class group of $K$,
 $Cl_{K,\ep}\cong I_{\ep}/P_{\ep,1}$ is the ray class group with
conductor $\ep$, and $\ov
 P_{\ep}=P_{\ep}/P_{\ep,1}$ is the principal ray with conductor $\ep$.
We have a short
 exact sequence
$$ 1 \longrightarrow \ov P_{\ep} \longrightarrow Cl_{K,\ep}
\longrightarrow Cl_K \longrightarrow 1. \tag*$$

Let $\ell$ be a prime dividing the orders of all three terms of (*), and
consider the exact sequence of $\ell$-primary components
$$ 1 \longrightarrow \ov P_{\ep}^{(\ell)} \longrightarrow
Cl_{K,\ep}^{(\ell)} \longrightarrow Cl_K^{(\ell)} \longrightarrow 1.
\tag*{${}_{\ell}$}$$
\smallskip
We are interested in primes $\ep$ for which the sequence (*${}_{\ell}$)
splits. Let $\ea_1,...,\ea_s \in I_{\ep}$ such that their images
$\ov\ea_i$ in $Cl_K^{(\ell)}$ form a basis of the finite abelian
$\ell$-group $Cl_K^{(\ell)}$.  Let $\ell^{m_i}$ be the order of
$\ov\ea_i$, $i=1,...,s$.  Then $\ea_i^{\ell^{m_i}}=(a_i)\in P_{\ep}$,
$i=1,...,s$.  Let $K_1=K(\root{\ell^{m_i}}\of{a_i}|1\leq i \leq
s)K(\mu_{\ell^m})$ with $m=\max\{m_1,...,m_s\}$.

\proclaim {Lemma 2.1} In order that the sequence
(*${}_{\ell}$) split, it is sufficient that $\ep$ split completely in
$K_1$. \endproclaim

\demo{Proof} Suppose $\ep$ splits in $K_1$.  Then $a_i$ is locally an
$\ell^{m_i}$-th power at $\ep$, so there exists $\alpha_i \in K^*$ such
that $\alpha_i^{\ell^{m_i}}\equiv a_i$ $(\mod \ep)$, $i=1,...,s$.  Set
$\eb_i:=\ea_i\alpha_i^{-1}$, $i=1,...,s$.  Then $\ov\eb_i=\ov\ea_i$,
$i=1,...,s$, and
$\eb_i^{\ell^{m_i}}=\ea_i^{\ell^{m_i}}\alpha_i^{-{\ell}^{m_i}}=(a_i)\alpha_i^{-{\ell}^{m_i}}$
with $a_i\alpha_i^{-\ell^{m_i}}\equiv 1$ (mod $\ep$), so $\eb_i^{\ell^{m_i}}\in
P_{\ep,1}$.  Let $\wt\eb_i$ be the image of $\eb_i$ in
$Cl_{K,\ep}^{(\ell)}$.  We have just seen that its order is
$\ell^{m_i}$.  Hence the $\wt\eb_i$ together generate a subgroup of
order at most $\ell^{\sum m_i}$.  Since the $\ov\eb_i=\ov\ea_i$ generated
the $\ell$-class group $Cl_K^{(\ell)}$, of order $\ell^{\sum m_i}$, it
follows that the $\eb_i$ together generate a complement to $\ov
P_{\ep}^{(\ell)}$ in $Cl_{K,\ep}^{(\ell)}$. \qed \enddemo

\bf Note. \rm In the proof of Lemma 2.1, we could have taken
$\ea_1,...,\ea_s \in I_{\ep}$
to be prime ideals outside of any given finite set of primes of $K$, by
virtue of the generalized Dirichlet density theorem.

 \head 3. $\ell^r$-torsion in the Brauer group
    \endhead

\proclaim{Theorem 3.1}  Let $K$ be a number field, $\ell$ a prime number,
$r$ a positive integer.  Then there exists an abelian $\ell$-extension
$L/K$ of exponent $\ell^r$ such that the local degree
$[LK_{\ep}:K_{\ep}]$ is equal to $\ell^r$ for every finite prime $\ep$ of
$K$, and is equal to $2$ at every real prime if $\ell=2$. \endproclaim

\demo{Proof}  Let $H_K$ be the Hilbert class field of $K$, and let
$H_K^{(\ell)}$ be the $\ell$-primary part of $H_K$.  Let $\ell^t$ be the
exponent of $Cl_K^{(\ell)}$. \smallskip Let $S$ be the set of primes of
$K$ that split completely in the extension $F$ of $K$ generated by
$H_K^{(\ell)}$, all $\ell^{r+t}$-th roots of all units (including roots
of unity) of $K$ , and the field $K_1$ of Lemma 2.1. For $\ep\in S$, $\ov
P_{\ep}\cong (\cO_K/\ep)^*$ modulo the image of the unit group $E_K$ of
$K$.  By definition of $S$, every unit of $K$ is an $\ell^{r+t}$-th power
in $K_{\ep}^*$, so the cyclic group $\ov P_{\ep}^{(\ell)}$ has order
divisible by $\ell^{r+t}$.  Let
 $R^{(\ep,\ell)}$ denote the $\ell$-ray class field with conductor
$\ep$. It follows that the subfield of $R^{(\ep,\ell)}$ fixed by the
$\ell^t$-torsion subgroup of $\Gal(R^{(\ep,\ell)}/K)$ is a cyclic
extension of $K$ of degree divisible by $\ell^r$.  Let $L^{\ep}$ be the
subfield of this field of degree exactly $\ell^r$ over $K$.  $L^{\ep}/K$
is totally ramified at $\ep$ and $\ep$ is the only prime of $K$ ramifying
in $L^{\ep}$.\smallskip

For each $\ep \in S$  the splitting of the short exact sequence
$(*_{\ell})$ allows us to fix a splitting map
$f_{\ep}:Cl_{K,\ep}^{(\ell)}\longrightarrow \ov P_{\ep}^{(\ell)}$ .
\smallskip

By \cite {\AT, Theorem 5, p. 105}, there exists a cyclic extension
$L_0/K$ of degree $\ell^r$ having local degree $\ell^r$ at the primes
$\el_1,...,\el_g$ dividing $\ell$, and local degree $2$ at the real
primes if $\ell=2$.
\smallskip
Let $\eq_1,...,\eq_e$ be the primes of $K$ not dividing $\ell$, which
ramify in $L_0$, ordered in such a way that $L_0/K$ has local degree
$\ell^r$ at $\eq_1,...,\eq_d$, $0\leq d \leq e$, and $L_0/K$ has local
degree less than $\ell^r$ at $\eq_{d+1},...,\eq_e$.  \smallskip Using the
Note following the proof of Lemma 2.1, we may assume without loss of
generality that the ideals $\ea_1,...,\ea_s$ in Lemma 2.1, representing
the elements of the class group, are prime ideals, distinct from
$\el_1,...,\el_g$ and from $\eq_1,...,\eq_e$.
\smallskip
Suppose $d<e$.  Let $\eq=\eq_{d+1}$, and let $\ell^a, \ell^b$ be the
ramification index and inertia degree respectively of $\eq$ in $L_0/K$.

We seek a prime $\ep\in S$  such that \smallskip (1)\ \ \ $\ep$ splits
completely in $L_0$, \smallskip (2)\ \ \ $\el_1,...,\el_g$ and
$\eq_1,...,\eq_d$ split completely in $L^{\ep}$, and \smallskip (3)\ \ \
the local degree of $L_0L^{\ep}/K$ at $\eq$ is $\ell^r$.  \smallskip (1)
is a Chebotarev splitting condition compatible with the Chebotarev
condition $\ep\in S$.  (2) is equivalent to a Chebotarev splitting
condition on $\ep$, by a ``reciprocity argument" (cf. \cite \KS),
formulated in the following lemma.

\proclaim {Lemma 2.2} Let $\ea$ be any prime of $K$, and let $\ep$ be any
prime in $S$ different from $\ea$.  Let $L'$ be a subfield of $L^{\ep}$
containing $K$, of degree $\ell^s$ over $K$.  Let $\ell^m$ be the order
of the image of $\ea$ in $Cl_K^{(\ell)}$, and let
$\ea^{\ell^m}=(\alpha)\in P_{\ep}$.  Then $\ea$ splits completely in $L'$
if and only if $\ep$ splits completely in
$K(\mu_{\ell^{m+s}},\root{\ell^{m+s}}\of{\alpha})$. \endproclaim

\demo{Proof}  Let $\widetilde{(\beta)}=f_{\ep}(\tilde\ea)\in \ov
P_{\ep}^{(\ell)}$, where $f_{\ep}:Cl_{K,\ep}^{(\ell)}\longrightarrow \ov
P_{\ep}^{(\ell)}$ is the splitting map above, and $\tilde\ea$ denotes the
image of $\ea$ in $Cl_{K,\ep}^{(\ell)}$. Then
$$\widetilde{(\beta)}^{\ell^m}={f_{\ep}(\tilde\ea)}^{\ell^m}=f_{\ep}(\tilde\ea^{\ell^m})
=f_{\ep}(\widetilde{(\alpha)})=\widetilde{(\alpha)}.$$  If  $\ep\in S$,
the cyclic group $ \ov P_{\ep}^{(\ell)}$ has order divisible by
$\ell^{r+t}$ hence by $\ell^{m+s}$ (because $m\leq t$ and $s\leq r$). Now
$\ea$ splits completely in $L'$ if and only if the Frobenius of $\ea$ in
$\Gal(L'/K)$ is trivial, which holds if and only if
$\widetilde{(\beta)}=f_{\ep}(\tilde\ea)$ is an $\ell^s$th power in $\ov
P_{\ep}^{(\ell)}$, which holds if and only if
$\widetilde{(\alpha)}=\widetilde{(\beta)}^{\ell^m}$ is an $\ell^{m+s}$th
power in $\ov P_{\ep}^{(\ell)}$ (since $ \ov P_{\ep}^{(\ell)}$ has order
divisible by $\ell^{m+s}$), which holds if and only if $\alpha$ is an
$\ell^{m+s}$th power mod $\ep$, which holds if and only if the polynomial
$x^{\ell^{m+s}}-\alpha$ has a root mod $\ep$, which holds if and only if
$x^{\ell^{m+s}}-\alpha$ factors into linear factors mod $\ep$   ($\ep$
splits completely in $K(\mu_{\ell^{r+t}})$), which holds if and only if
$\ep$ splits completely in
$K(\mu_{\ell^{m+s}},\root{\ell^{m+s}}\of{\alpha})$. \qed \enddemo It
follows from the Lemma that the second condition is a Chebotarev
splitting condition on $\ep$, hence compatible with the previous
conditions.  The third condition, that  the local degree of
$L_0L^{\ep}/K$ at $\eq$ is $\ell^r$, will hold if $\eq$ is unramified in
$L^{\ep}/K$ with degree of inertia exactly $\ell^{r-a}$.  This condition
is equivalent, by Lemma 2.2, to the condition that $\ep$ splits
completely in $K(\mu_{\ell^{m+a}},\root{\ell^{m+a}}\of{\alpha})$ but not
in $K(\mu_{\ell^{m+a+1}},\root{\ell^{m+a+1}}\of{\alpha})$, where $\eq$
plays the role of $\ea$ in Lemma 2.2.  $m$ is the order of $\eq$ in the
class group.  Note that $a<r$ here, so $m+a+1 \leq r+t$.  In order that
this last Chebotarev condition be compatible with the preceding ones, it
is necessary and sufficient that $\root{\ell^{m+a+1}}\of{\alpha}$ not lie
in the composite of all the fields in which $\ep$ has been required to
split.  This follows from ramification considerations, since the
ramification index of $\eq$ in the composite of all the fields in which
$\ep$ has been required to split, is exactly $\ell^a$, whereas the
ramification index of $\eq$ in
$K(\mu_{\ell^{m+a+1}},\root{\ell^{m+a+1}}\of{\alpha})$ is $\ell^{a+1}$.
We have therefore proved the existence of a prime (in fact infinitely
many) $\ep=\ep_1$ satisfying the three conditions. \smallskip Similarly,
there exists a prime $\ep_2\in S$, such that  \smallskip (4)\ \ \ $\ep_2$
splits completely in $L_0$ and $L_1:=L^{\ep_1}$,
\smallskip (5)\ \ \ $\el_1,...,\el_g,\eq_1,...,\eq_d,\eq_{d+1},\ep_1$ split
completely in $L_2:=L^{\ep_2}$, and
\smallskip (6)\ \ \ the local degree of $L_0L^{\ep_2}/K$ at $\eq_{d+2}$ is
$\ell^r$. \bigskip  Proceeding in this manner, we get
$\ep_1,...,\ep_{e-d}$, such that
\smallskip (7)\ \ \ $\ep_i$ splits completely in $L_j$ for $i\neq j$, \smallskip
(8)\ \ \ $\el_1,...,\el_g,\eq_1,...,\eq_e$ split completely in
$L_1,...,L_{e-d}$, and \smallskip (9)\ \ \ the local degree of the
composite $L_0L_1\cdots L_{e-d}$ is equal to $\ell^r$ at each of the
primes $\el_1,...,\el_g,\eq_1,...,\eq_e,\ep_1,...,\ep_{e-d}$.
\smallskip
Let us now extend $\el_1,...,\el_g,\eq_1,...,\eq_e,\ep_1,...,\ep_{e-d}$
to an enumeration

$$\el_1,...,\el_g,\eq_1,...,\eq_e,\ep_1,...,\ep_{e-d},\ea_1,\ea_2,...,\ea_s,\ea_{s+1},...$$
of all the finite primes of $K$, where $\ea_1,\ea_2,...,\ea_s$ are the
primes of Lemma 2.1.  \smallskip We now seek a prime $\ep=\ep_{e-d+1}\in
S$ such that \smallskip (10)\ \ \ $\ep$ splits completely in
$L_0,L_1,...,L_{e-d}$,
\smallskip (11)\ \ \ $\el_1,...,\el_g,\eq_1,...,\eq_e,\ep_1,...,\ep_{e-d}$ split
completely in $L^{\ep}$, and  \smallskip (12)\ \ \ $\ea=\ea_1$ is inert
in $L^{\ep}$.
\smallskip (10) and (11) are equivalent to saying
that  $\ep$ splits completely in a Galois extension (a composite of
$L_i$'s and fields used in proving compatibility of previous Chebotarev
conditions) in which $\ea$ is unramified.
 (12) is equivalent to the condition that $\ea$ does
not split completely in the subfield of $L^{\ep}$ of degree $\ell$ over
$K$. By Lemma 2.2 applied to this $\ea$, this condition is equivalent to
the (Chebotarev) condition that $\ep$ does not split completely in
$K(\mu_{\ell^{m+1}},\root{\ell^{m+1}}\of{\alpha})$, where $m$ corresponds
to this $\ea$.  Since $\ea$ ramifies in
$K(\mu_{\ell^{m+1}},\root{\ell^{m+1}}\of{\alpha})$, this last Chebotarev
condition is compatible with the others.  \smallskip  We remove
$\ep=\ep_{e-d+1}$ from the enumeration, denote $L^{\ep}$ by $L_{e-d+1}$,
and move to the next ideal $\ea_2$. A similar argument yields a prime
$\ep=\ep_{e-d+2}$ and a corresponding $L_{e-d+2}$, such that
\smallskip (13)\ \ \ $\ep$ splits completely in $L_0,L_1,...,L_{e-d+1}$,
\smallskip (14)\ \ \ $\el_1,...,\el_g,\eq_1,...,\eq_e,\ep_1,...,\ep_{e-d+1}$ split
completely in $L^{\ep}$, and \smallskip (15)\ \ \ $\ea=\ea_2$ is inert in
$L^{\ep}=L_{e-d+2}$.  \smallskip Continuing in this way, a sequence
$\{\ep_i\}$ of primes of $K$ is generated, with a corresponding sequence
of fields $\{L_i\}$.

Form the composite $L:=L_0L_1L_2\cdots$. The local degree of $L/K$ at
every finite prime is $\ell^r$. In fact, the local degree of
$L_0L_1\cdots L_{e-d}$ is equal to $\ell^r$ at each of the primes
$\el_1,...,\el_g,\eq_1,...,\eq_e,\ep_1,...,\ep_{e-d}$, as we showed
earlier; these primes split completely in $L_j$ for $j>e-d$; at every
$\ep_i$, the local extension is cyclic and totally ramified of degree
$\ell^r$; and at every $\ea\notin
\{\el_1,...,\el_g,\eq_1,...,\eq_e,\ep_1,\ep_2,...\}$, the local extension
is cyclic and unramified of degree $\ell^r$. The local degree at the real
primes is $2$ since this is true already for $L_0$. \qed \enddemo

\head 4. The function field case. \endhead

Let $K$ be a global function field of characteristic $p$.  Let $\ell$ be
a prime different from $p$. Theorem 3.1 holds also for
$K$, with essentially the same proof, except simpler, since there are no
primes dividing $\ell$.  We indicate how this works. \smallskip

By Chebotarev's density theorem for function fields (see e.g.
\cite\Weil), there exists a prime $\ep$ which is inert in the constant
degree $\ell$ extension of $K$. This property is equivalent to $\ep$
having degree prime to $\ell$ (see e.g. \cite\Po). Fix one such prime
$\ep_{\infty}$. Let $\cO$ be the ring of all elements of $K$ which are
integral at all primes $\ep\neq \ep_{\infty}$ of $K$. $\cO$ is a
Dedekind
domain.  Call $\ep_{\infty}$ the infinite prime of $K$ and all the
others
finite primes of $K$.
There is an exact sequence $$1\longrightarrow \Pic^0(K) \longrightarrow
\Pic(\cO) \longrightarrow \dZ/d_{\infty}\dZ \longrightarrow 0$$ with
$\Pic^0(K)$ the group of divisor classes of degree zero of $K$,
$\Pic(\cO)$ the ideal class group of $\cO$, and the third arrow is the
degree mod $d_{\infty}$ map, where $d_{\infty}$ denotes the degree of
$\ep_{\infty}$.  As $d_{\infty}$ and $\ell$ are relatively prime, the
$\ell$-primary components of $\Pic^0(K)$ and $\Pic(\cO)$ are isomorphic;
we denote them by $Cl_K^{(\ell)}$.  If $\ep$ is a finite prime of $K$,
let $R^{(\ep)}$ be the maximal abelian extension of $K$ with conductor
dividing $\ep$, in which $\ep_{\infty}$ splits completely, i.e. the ray
class field mod $\ep$. $\Gal(R^{(\ep)}/K)$  is canonically isomorphic,
via the reciprocity map, to $Cl_{K,\ep}$, the  ray class group mod $\ep$
(See \cite {\hayes, p. 204}.

The proofs of Lemma 2.1, Lemma 2.2 and Theorem 3.1 proceed as in the
number field case, except that there are no primes
$\el_1,...,\el_g,\eq_1,...,\eq_e$ to deal with; the field $L_0$ is
replaced with  the constant degree $\ell^r$ extension of $K$, and $e=0$.
$L$ is the composite of the fields $L_0,L_1,L_2,...$.  We therefore have

\proclaim {Theorem 4.1} Let  $K$ be a global function field of
characteristic $p$, $\ell$  a prime different from $p$, $r$ a positive
integer.   Then there exists an abelian $\ell$-extension $L/K$ of
exponent $\ell^r$ such that the local degree $[LK_{\ep}:K_{\ep}]$ is
equal to $\ell^r$ for every prime $\ep$ of $K$. \endproclaim

 \head 5. The $n$-torsion subgroup of the Brauer
group of $K$ \endhead

\proclaim {Theorem 5.1} Given a global field $K$  and a positive integer
$n$,  there exists an abelian extension $L/K$ (of exponent $n$) such that
the $n$-torsion subgroup of the Brauer group of $K$ is equal to the
relative Brauer group of $L/K$.

\endproclaim

\demo{Proof} Consider the case $n=\ell^r$, $\ell$ prime.  By Theorems 3.1
and 4.1, if $\charr(K)\neq \ell$, there exists an abelian
$\ell$-extension $L/K$ whose local degree at every nonarchimedian prime
is equal to $\ell^r$, and is equal to $2$ at the real primes if $K$ is a
number field and $\ell=2$. If $\charr(K)=\ell$, the same result holds by
\cite{\Po}. It follows from the fundamental theorem of class field theory
on the Brauer group of a number field that $L$ splits every algebra class
of order dividing $\ell^r$, and conversely, any algebra class split by
$L$ has order dividing $\ell^r$. For general $n$, the theorem follows
from a straightforward reduction to the prime power case (see
\cite{\as}). \qed\enddemo
\smallskip
\bf Remark. \rm Adrian Wadsworth has pointed out to us that in the case
$\charr(K)=\ell$, $Br_{\ell^r}(K)$ is the relative Brauer group of a
purely inseparable extension, namely $L=\root{\ell^r}\of{K}$, by results
of Albert \cite {\albert, p. 109}. \smallskip
\bf Acknowledgment. \rm We are grateful to Harold Stark for valuable
discussions.

\end